\documentclass[12pt]{amsart}

\newtheorem{theorem}{Theorem}[section]

\newtheorem{proposition}[theorem]{Proposition}

\theoremstyle{definition}

\newtheorem{remark}[theorem]{Remark}

\theoremstyle{parrafo}

\numberwithin{equation}{theorem}

\begin{document}

\title[]{Dimension dependency of
 the weak type $(1,1)$ bounds for  maximal functions associated to
finite radial measures}

\author{J. M. Aldaz}
\address{Departamento de Matem\'aticas y Computaci\'on,
Universidad  de La Rioja, 26004 Logro\~no, La Rioja, Spain.}


\thanks{Partially supported by Grant BFM2003-06335-C03-03 of the D.G.I. of
Spain}







\begin{abstract} We show that the lowest constants appearing
in the weak type (1,1) inequalities
satisfied by the centered Hardy-Littlewood maximal function
associated to  some finite radial measures, such as the standard
gaussian measure, grow exponentially fast with the dimension.
\end{abstract}


\maketitle


\section {Introduction}

\markboth{J. M. Aldaz}{Weak type $(1,1)$ bounds}

Let $M$ be the centered maximal operator (cf. (\ref{HLMF}) below) associated to euclidean
balls and Lebesgue measure. It is well known that if $1 < p \le \infty$, then there exists
a constant $c_p$ such that for all $f\in L^p (\mathbb R^d)$, we have $\|Mf\|_p\le c_p \|f\|_p$.
When $p =\infty$, trivially $c_p = 1$. The standard proof of $\|Mf\|_p\le c_p \|f\|_p$
for $1<p<\infty$, via weak type (1,1) inequalities and interpolation, rests on covering lemmas
of a geometric character. This fact leads to values of $c_p$ that grow exponentially
with the dimension $d$. An alternative proof using the method of rotations gives constants 
$c_p$ whose growth is linear in $d$ (but yields no weak type (1,1) inequality). 
In this context,  E. M. Stein was able to show that in fact one can take $c_p$ to be independent of $d$ (\cite{St1}, \cite{St2}, see also \cite{St3}).
 A motivation for the study of  $L^p$ bounds uniform in $d$ comes from the desire to
extend (at least some parts of) harmonic analysis in $\mathbb
R^d$, to the infinite dimensional case. 
Stein's result was
generalized to the maximal function defined using an arbitrary norm by J. Bourgain (\cite{Bou1},
\cite{Bou2}, \cite{Bou3}) and A. Carbery (\cite{Ca}) when $p>3/2$.
For $\ell_q$ balls, $1\le q <\infty$, D. M\"{u}ller \cite{Mu}
showed that uniform bounds again hold for every $p
> 1$ (given $1\le q <\infty$, the $\ell_q$ balls
are defined using the norm $\|x\|_q :=\left( x_1^q+ x_2^q+\dots + x_d^q\right)^{1/q}$). With respect to weak type bounds, in \cite{StSt} E. M. Stein
and J. O. Str\"{o}mberg proved, among other things, that the smallest (i.e., the
best)
constants in the weak type (1,1) inequality satisfied by $M$ grow
at most like $O(d)$, and asked if uniform bounds could be found.
Since then, there has been remarkably little progress on this
question (see, for instance, \cite{A}, \cite{AV} for the case of
cubes).

Here we study the weak type (1,1) problem for integrable radial
densities defined via bounded decreasing functions, the canonical
example being the standard gaussian measure. This is a natural
variant of Stein and Str\"{o}mberg's question, given the growing
interest in what has been termed ``gaussian harmonic analysis",
where Lebesgue measure is replaced by  the standard gaussian
measure,  and also because of the importance of gaussian measures
and other probabilities in the infinite dimensional setting (see,
for instance, \cite{Bo}). For the measures considered in this
paper, instead of uniform bounds we have exponential increase: if
$\mu$ is a finite radial Borel measure on $\mathbb R^d$ defined
by a bounded decreasing function $f$, and if $c_d$ denotes the
smallest constant appearing in the weak type (1,1) inequality
satisfied by the associated maximal function $M_{\mu}$, then for every
$d$ we have $c_d \ge \left(1 + \frac{ 2 \sqrt 2}{\sqrt{ 3\pi d}}\sqrt{1 + \frac1d}\right)^{-1} \left(
\frac{2}{\sqrt 3}\right)^{d/6}.$

\section
{Notation and results}

Given a locally finite Borel measure $\mu$ on $\mathbb R^d$ (so compact sets have
finite measure) and a locally integrable function $f$, 
the associated
centered maximal function $M_{\mu}f$ is defined by
\begin{equation}\label{HLMF}
M_{\mu} f(x) := \sup _{\{r > 0: \mu (B(x, r)) > 0\}} \frac{1}{\mu
(B(x, r))} \int _{B(x, r)} \vert f\vert d\mu,
\end{equation}
where $B(x, r)$ denotes the euclidean {\it closed} ball of radius
$r > 0$ centered at $x$ (the choice of closed balls in the
definition is mere convenience; using open balls instead does not
change the value of $ M_{\mu} f (x)$). The boundary of $B(x, r)$
is the sphere $\mathbb S (x,r)$. Sometimes we use $B^d (x, r)$ and
$\mathbb S^{d-1} (x,r)$ to make their dimensions explicit. If
$x=0$ and $r=1$, we just write $B^d$ and $\mathbb S^{d-1}$.

It is a consequence of Besicovitch covering theorem that there
exists a constant $c= c(d)$, independent of $\mu$, such that for
every $f \in L^1(\mathbb R^d, \mu )$ and every $\alpha > 0$, we
have $\alpha \mu (\{ M_{\mu} f \ge \alpha \}) \le c \Vert f\Vert
_1$. In fact, by the Theorem in pg. 227 of \cite{Su} we may take
$c = (2.641 + o(1))^d$. Of course, if we are interested in the
lowest such $c$, then $c= c(\mu )$ will depend on $\mu$.
Note that it makes no difference in the determination of the smallest
constant if instead of the strict inequality $\{ M_{\mu} f >
\alpha \}$ we use $\{ M_{\mu} f \ge \alpha \}$.

The following easy result highlights the fact that uniform bounds will often fail to
exist when dealing with sequences $\{\mu_d\}_{d= 1}^\infty$ of  measures satisfying 
$\mu_d(\mathbb R^d) < \infty$, due to the
decay this condition imposes. We
take $\frac{a}{\mu (B(y, r))}$ to mean $\infty$  when
$a > 0$ and $\mu (B(y, r)) =0$.

\begin{proposition}\label{proposition} Let $\mu$ be a locally finite Borel measure  on $\mathbb
R^d$, and let $c_d$ be the smallest constant appearing in the weak
type (1,1) inequality satisfied by $M_{\mu}$. Given any ball 
$B (x,r)$ with $x$ in the support of $\mu$,
$$
c_d \ge \inf_{y\in \mathbb S (x,r)} \frac{\mu( B(x, r))}{\mu (B(y,
r))}.
$$

\end{proposition}

\noindent{\em Proof.} Fix $r > 0$ and let $x$ belong to the support of $\mu$.
By a standard approximation argument we may consider, instead of a
function, the Dirac delta $\delta_x$ placed at $x$. Thus $c_d \ge
\sup_{\alpha > 0} \alpha \mu (\{M_{\mu} \delta_x \ge \alpha \}).$
Note that for some $y\in \mathbb S (x,r)$, $\mu (B(y, r)) > 0$. 
Also, if $y\in \mathbb S (x,r)$ and $\mu (B(y, r)) = 0$, then
$M\delta_x (y) = \infty$, since for all $s > r$,  $\mu (B(y, s)) > 0$ and
$\lim_{s\downarrow r} \mu (B(y, s)) = 0$.  Let $\alpha_0 := \inf_{y\in \mathbb S
(x,r)} \frac{1}{\mu (B(y, r))}$. Then
$M_{\mu} \delta_x (z)\ge \alpha_0$ for every  $z = (1-t) x + t y$, $y\in \mathbb S
(x,r)$, $0 < t \le 1$,  since $x \in B(z,
tr)\subset B(y, r)$. It follows that 
$B(x,r) \subset \{M_{\mu} \delta_x \ge \alpha_0 \}$, so $c_d \ge
\alpha_0 \mu ( B(x,r))$, as claimed. \qed

\begin{remark}\label{trivial} Let $\mu$ be a locally finite Borel measure on $\mathbb{R}^d$.
By the Lebesgue Theorem on differentiation of integrals, 
$M_\mu f(x)\ge |f|(x)$ for $\mu$ a.e. $x$. Now, fix $\varepsilon > 0$, and let $f :=(1 +\varepsilon) \chi_{B(0,1)}$. Then $\mu (B(0,1)) =\mu(\{f>1\}) \le (1+\varepsilon) c_d\int \chi_{B(0,1)}d\mu$,
so $c_d\ge 1$.
\end{remark}

Fix $d\in \mathbb{N}\setminus 0$. 
Let $f: [0,\infty ) \to [0,\infty )$ be a nonincreasing
 function,  let $\sigma^{d-1}$ denote the
area on the unit sphere $\mathbb S^{d-1}$, let $\tilde
\sigma^{d-1}$ denote the {\em normalized} area on $\mathbb
S^{d-1}$ (thus $\tilde \sigma^{d-1}$ is a probability), and let
$\lambda ^d$ be Lebesgue measure on $\mathbb R^d$. Then the function $f$
defines a  rotationally invariant (or radial) measure
$\mu$ via
\begin{equation}\label{defrad}
\mu (A) := \int_A  f(|y|)  d\lambda^{d} (y).
\end{equation}
We remark that since $f$ is nonincreasing, it is bounded by $f(0)$. Additionally, we
shall assume that $f$ is not $0$ a.e., so $\mu (\mathbb{R}^d)>0$, and furthermore,
that $f(x) x^{d-1}\in L^1 ([0,\infty))$, so $\mu (\mathbb{R}^d)<\infty$, as can be seen by
integrating in polar coordinates.

\begin{theorem}\label{theorem} Fix $d\in \mathbb{N}\setminus 0$.
 Let $f: [0,\infty ) \to [0,\infty )$ be a nonincreasing
 function and let $\mu$ be the radial measure defined via (\ref{defrad}). Assume  $f$ is
 such that  $0 < \mu (\mathbb{R}^d)<\infty$, and let
$c_d$ be the smallest constant appearing in the weak type (1,1)
inequality satisfied by $M_{\mu}$. Then 
\begin{equation}\label{best}
c_d \ge \left(1 + \frac{ 2 \sqrt 2}{\sqrt{ 3\pi d}}\sqrt{1 + \frac1d}\right)^{-1} \left(
\frac{2}{\sqrt 3}\right)^{d/6}.
\end{equation}

\end{theorem}

\noindent{\em Proof.}  Since $c_1 \ge 1$ by Remark \ref{trivial}, the lower bound (\ref{best}) holds for $d=1$, so we may assume that
$d\ge 2$. Given a unit vector $v\in \mathbb
R^{d}$ and $\varepsilon \in [0, 1)$, the $\varepsilon$ spherical
cap about $v$ is the set $C(\varepsilon , v) :=\{\theta \in
\mathbb S^{d-1}: \langle \theta, v\rangle \ge \varepsilon\}$. Note
that spherical caps are just geodesic balls $B_{\mathbb
S^{d-1}}(x, r)$ in $\mathbb S^{d-1}$. In the special case $v=e_1 =
(1,0,\dots , 0)$, $\varepsilon =2^{-1}$, we have $C(2^{-1}, e_1) =
B_{\mathbb S^{d-1}}(e_1, \pi/ 3)$. Next we remind the reader of
some well-known facts that will be used in the sequel: i)
$\lambda^d (B^d) =\frac{\pi^{d/2}}{\Gamma (1 + d/2)}$; ii)
$\sigma^{d-1} (\mathbb S^{d-1}) = d \lambda^d (B^d)$; iii)
 $\sigma^{d-1} (B_{\mathbb
S^{d-1}} (x,r)) = \sigma^{d-2} (\mathbb S^{d-2})\int_0^r
\sin^{d-2}t dt$ (see, for instance, (A.11) pg. 259 of \cite{Gra}
for a more general statement). We shall also use the following known 
 and elementary 
estimate (cf. Exercise 5, pg. 216 of \cite{Web}). 
The short derivation is included for the reader's convenience. Recall that on  $\{x > 0\}$,
 $\Gamma$ is log-convex  and  $\Gamma (x + 1) = x \Gamma (x)$. Writing
$d+2 = 2^{-1} (d+1) + 2^{-1} (d+3)$, we have
\begin{equation}\label{ratio}
\frac{\Gamma (1 + d/2)}{\Gamma (1/2 + d/2)}
\le \frac{\left(\Gamma (2^{-1} (d+1))\right)^{1/2}\left(\Gamma (2^{-1} (d+3))\right)^{1/2}}{  \Gamma (2^{-1} (d+1))} = \left(\frac{d+1}{2}\right)^{1/2}.
\end{equation}

Note that from i), ii), iii), (\ref{ratio})
and the fact that $\cos t \ge 1/2$ on $[0,\pi/3]$, we get the
following upper bound on the normalized area of $C(2^{-1}, e_1)$:
\begin{equation}\label{in}
\tilde \sigma^{d-1}\left( C(2^{-1}, e_1)\right) \le 2 \frac{\sigma^{d-2}
(\mathbb S^{d-2})}{\sigma^{d-1}(\mathbb S^{d-1})} \int_0^{\pi/3}
\sin^{d-2}t \cos t dt
\end{equation}
\begin{equation}\label{in2}
= \frac{2}{d} \frac{\lambda^{d-1} (B^{d-1})}{\lambda^{d} (B^{d})}
\left( \frac {\sqrt 3}{2}\right)^{d-1}  \le
 \left( \frac
{\sqrt 3}{2}\right)^{d}  \frac{ 2 \sqrt 2}{\sqrt{ 3\pi d}}\sqrt{1 + \frac1d}.
\end{equation}

Note next that the function $h(R):= \frac{\mu (B(0, R))}{\mu (B(0, (\sqrt 3 /2)
R))}$ is continuous, and $\lim_{R\to\infty} h(R) = 1$ by the
finiteness of the measure. It follows that there is a largest real
number $R_1$ such that $h(R_1) = (2/\sqrt 3 )^{d/6}$, provided of
course that the set $\{h \ge (2/\sqrt 3 )^{d/6}\}$ is nonempty. To see that this is always the case,  note that since $f$ is nonincreasing, the same happens with the
averages $\frac{1}{\lambda^d (B(0,R))} \int_{B(0,R)} f(|x|) dx$. Thus 
$\lim_{R\to 0} \frac{1}{\lambda^d (B(0,R))} \int_{B(0,R)} f(|x|) dx = L$ exists and
$L\le f(0)< \infty$. It
follows that 
$$\lim_{R\to 0} h(R) =  
\lim_{R\to 0} \frac{\frac{\lambda^d (B(0,R))}{\lambda^d (B(0,R))} \int_{B(0,R)} f(|x|) dx}{\frac{\lambda^d (B(0,(\sqrt 3 /2)R))}{\lambda^d (B(0,(\sqrt 3 /2)R))} \int_{B(0,(\sqrt 3 /2)R)} f(|x|) dx}=
  (2/\sqrt 3 )^{d}.
$$

By rotational invariance and the previous proposition, in order to prove the
theorem
it is enough to check that  
$$\frac{\mu (B(0, R_1))}{\mu (B(R_1 e_1,  R_1))}
\ge \left(1 + \frac{ 2\sqrt 2}{\sqrt{ 3\pi d}}\sqrt{1 + \frac1d}\right)^{-1}
\left( \frac{2}{\sqrt 3}\right)^{d/6}.
$$
We split $\mu (B(R_1 e_1,  R_1)) = \mu (B(0, R_1)\cap B(R_1 e_1,
R_1))+ \mu ((B(0, R_1)^c\cap B(R_1 e_1,  R_1))$ and estimate
each of the summands.  Note that
$$
B(0, R_1)\cap B(R_1 e_1, R_1) \subset B(2^{-1} R_1 e_1, 2^{-1}
\sqrt 3 R_1).
$$
For every pair of points $(x,y)$ with $x\in B(0, R_1)\setminus
B(R_1 e_1, R_1)$, $y\in B(R_1 e_1, R_1)\setminus B(0, R_1)$, we
have  $\vert x \vert < \vert y \vert$, so $f(\vert x \vert) \ge
f( \vert y \vert)$. By the choice of $R_1$ and the preceding
observation,
\begin{equation}\label{inside}
\mu (B(0, R_1)\cap B(R_1 e_1, R_1)) \le \mu (B(2^{-1} R_1 e_1,
2^{-1} \sqrt 3 R_1)) 
\end{equation}
$$
 \le \mu ( B(0, 2^{-1} \sqrt 3 R_1))
= (\sqrt 3/ 2)^{d/6} \mu (B(0, R_1)).
$$
Let $E$ be the semi-cone in $\mathbb R^d$ defined by $x_1 =
3^{-1/2}\sqrt{x_2^2 + \dots + x^2_d}$, and let $E^\prime := \{x_1
\ge 3^{-1/2}\sqrt{x_2^2 + \dots + x^2_d}\}$ denote the solid semi-cone
determined by $E$. Then
$$B(0, R_1)^c\cap B(R_1 e_1,
R_1)\subset B(0,(2/ \sqrt 3)^5 R_1) \cap E^\prime.$$
  By
rotational invariance of the probability measure
$$\nu (A) :=
\frac{\mu (B(0,(2/ \sqrt 3)^5 R_1) \cap A)}{\mu (B(0,(2/ \sqrt
3)^5 R_1))},$$
 the $\nu$-measure of  $E^\prime$ is just the normalized
 area of its intersection with
 the sphere, i.e.
$$
\nu (E^\prime) =  \tilde \sigma^{d-1} (E^\prime\cap \mathbb
S^{d-1}) = \tilde \sigma^{d-1}\left( C(2^{-1}, e_1)\right).
$$
From  the choice of $R_1$  and the upper bound (\ref{in}, \ref{in2}) on $\tilde
\sigma^{d-1}\left( C(2^{-1}, e_1)\right)$,  we get
\begin{equation}\label{out}
\mu ((B(0, R_1)^c\cap B(R_1 e_1,  R_1)) \le  \nu (E^\prime)
\mu (B(0,(2/ \sqrt 3)^5 R_1))
\end{equation}
\begin{equation*}
< \tilde \sigma^{d-1}\left( C(2^{-1}, e_1) \right)(2/\sqrt 3 )^{5d/6} \mu
(B(0, R_1))
\end{equation*}
\begin{equation*}
\le  \left( \frac {2}{\sqrt 3}\right)^{5d/6} \mu (B(0, R_1))
\left( \frac {\sqrt 3}{2}\right)^{d} \frac{ 2 \sqrt 2}{\sqrt{3
\pi d}}\sqrt{1 + \frac1d}
\end{equation*}
\begin{equation*}
= \mu (B(0, R_1)) \left( \frac {\sqrt 3}{2}\right)^{d/6}
\frac{ 2\sqrt 2}{\sqrt{ 3\pi d}}\sqrt{1 + \frac1d}.
\end{equation*}
Putting   together the estimates starting at (\ref{inside}) and at (\ref{out}), we obtain
$$
\frac{\mu (B(0, R_1))}{\mu (B(R_1 e_1,  R_1))} \ge \frac{\mu
(B(0, R_1))}{\left(1 + \frac{ 2\sqrt 2}{\sqrt{ 3\pi d}}\sqrt{1 + \frac1d}\right)
\left( \frac{\sqrt 3}{2}\right)^{d/6} \mu_d B(0, R_1)},
$$
as desired. \qed

\begin{remark} For each $d\in \mathbb{N}\setminus 0$ let $f_d: [0,\infty ) \to [0,\infty )$ be a nonincreasing
integrable function, not $0$  a.e., such that $f_d (x) x^{d-1}$ is integrable, 
and let $\mu_d$ be the nontrivial finite radial measure on $\mathbb{R}^d$ defined using $f_d$. By the preceding theorem, there cannot be a uniform bound for the lowest constants appearing in the
weak type (1,1) inequalities satisfied by  $\mu_d$. Note that for different values
of $d$ the functions $f_d$ might be totally unrelated, or on the other extreme, might always
be the same function $f$ (provided $f (x) x^{d-1}$ is integrable for all $d$). But this has no effect on the lack of uniform bounds.
\end{remark}

\begin{remark}\label{oned} When $d = 1$, one can easily improve on the trivial bound $c_1 \ge 1$. Fix $\varepsilon > 0$ and choose $R>0$ such that
$\frac{\mu ([0,2 R])}{\mu ([0, R])} < 1 + \varepsilon$. Then 
$\frac{\mu ([-R, R])}{\mu ([0, 2 R])} = \frac{\mu ([-R, R])}{\mu ([-2 R, 0])} > 
\frac{2 \mu ([0, R])}{(1+\varepsilon)\mu ([0,  R])}=\frac{2}{1 + \varepsilon}$. So in fact $c_1 \ge 2$ for all finite nontrivial radial measures on $\mathbb{R}$.
\end{remark}

\begin{remark}  We briefly comment on the role that different assumptions
play in the  proof of Theorem \ref{theorem}. Finiteness  of $\mu$ was only
used to show that there is an $R
> 0$ such that $h(R)\ge (2/\sqrt 3 )^{d/6}$ and $h((2/\sqrt 3
)^{i} R)\le (2/\sqrt 3 )^{d/6}$ for $i=1,2,3,4,5$. So exponential dependency holds
for
 families of radial measures, not
necessarily finite, provided there exists a ball in each dimension with center in the
support of the corresponding measure, and a similar kind of decay (clearly the proof can be adapted to
decays lower than the one considered above).  Normalization of the
measures is not an issue either: if $c> 0$, then $M_{\mu} = M_{c
\mu}$. The radial assumption makes it easy to check the decay and
to apply Proposition \ref{proposition}, but as noted above, the
existence of some ball with such decay is likely to be prevalent
even among nonradial (finite)
 measures. Regarding the assumption of absolute continuity of the measures $\mu_d$ in $\mathbb{R}^d$, it
 is a simple way to ensure that in some sense the dimension
of $\mu_d$ goes to infinity as $d\to\infty$ 
(cf. for instance Chapter 10 of \cite{Fal} for a definition of dimension of a measure, and in particular, Prop. 10.5 pg. 174: it follows that 
if $\mu_d$ is absolutely continuous in $\mathbb{R}^d$ then 
$\dim \mu_d = d$).   A simple example where $\dim \mu_d$ does not grow with
$d$ is obtained by setting   $\mu _d = \delta_0$ for every
$d$. Here  we always have $\dim \mu_d = 0$ and $c_d = 1$. 
It is also easy to  obtain
 singular examples where there is  dependency of
$c_d$ on $d$: denote respectively by $\gamma_d$ and by $\delta_0$ the standard gaussian measure on $\mathbb R^d$ and the point mass at the origin of $\mathbb R^d$. Then $\gamma_d\times \delta_0$ on $\mathbb
R^{2d}$ is singular, $d$-dimensional, and $c_d$ grows exponentially with $d$.

  Finally, we point out that there are some
natural families of singular radial finite measures to which the
 arguments in the proof of Theorem \ref{theorem} do not apply, such as, for instance, area  on
$\mathbb S^{d-1}$.
\end{remark}

\begin{remark} When dealing with  concrete families of measures, the
additional information may lead  to  more precise bounds. For
instance, let $\nu_d (A) := \lambda^d (A\cap B^d)$ be Lebesgue
measure restricted to the unit ball in $\mathbb R^d$. By Remark (\ref{oned}), $c_1(\nu_1) \ge 2$, and for $d\ge 2$, the argument used in
(\ref{inside}) immediately gives $c_d (\nu_d)
\ge (2/ \sqrt 3)^d$. Nevertheless, it is possible to do   better by
estimating directly the volume of the solid cap $B^d \cap \{x_1 \ge
2^{-1}\}$:
\begin{equation*}
 \nu_{d} (B( e_1, 1)) = 2 \lambda^d (B^d \cap \{x_1 \ge
2^{-1}\}) = 2 \lambda^{d-1} (B^{d-1})\int_{1/2}^1\left(\sqrt{1-x_1^2}\right)^{d-1}dx_1
\end{equation*}
\begin{equation*}
 \le 4 \lambda^{d-1} (B^{d-1}) \int_{\pi/6}^{\pi/2}\cos^d t \sin t dt =   \frac{ 4}{ d + 1}\left( \frac{ \sqrt 3}{2}\right)^{d+1} \lambda^{d-1} (B^{d-1}).
 \end{equation*}
Using Proposition \ref{proposition} and (\ref{ratio}) we have 
\begin{equation}\label{stein}
 c_d (\nu_d) \ge \frac{\lambda^{d} (B^{d})}{\nu_{d} (B( e_1, 1))}\ge
\frac{\sqrt{\pi (d+1)}}{\sqrt 6} \left( \frac
{2}{\sqrt 3}\right)^{d}.
\end{equation}
The same estimate holds if $R > 0$ is chosen arbitrarily and we define $\nu_d (A) := \lambda^d (A\cap B(0, R))$,  since 
$\frac{\lambda^{d} (B(0,R)}{\nu_{d} (B( R e_1, R))}$ is independent of $R$.
Now fix $d$. As noted in the previous remark, multiplying $\nu_d$ by a nonzero constant does not
change any estimate; set  $\nu_n (A) := C_n \lambda^d (A\cap B(0, 1/n))$,
where $C_n$ is chosen to make $\nu_n$  a probability. Then $\nu_n\to \delta_0$ in the weak*
topology (weakly, in the standard probabilistic terminology), but
 by Remark \ref{oned}, the lower bound (\ref{stein}) and the preceding comments (or by Theorem \ref{theorem} if $d$ is high
 enough), $\liminf_n c_d(\nu_n) > 1= c_d(\delta_0 )$.
A similar observation can be made about the sequence defined by
$\nu_n (A) := \lambda^d (A\cap B(0, n))$  and  $\lambda^d$.
\end{remark}

\end{document}